\documentclass[12pt]{article}
\usepackage{amsmath, amssymb, amsfonts, amsthm}
\def\CC{{\rm \kern.24em \vrule width.02em height1.4ex
depth-.05ex \kern-.26em C}}

\def\TagOnRight

\def\AA{{\it I}\hskip-3pt{\tt A}}

\def\QQ{\rlap {\raise 0.4ex \hbox{$\scriptscriptstyle |$}}
  {\hskip -0.1em Q}}

\newcommand{\be}{\begin{equation}}
\newcommand{\ee}{\end{equation}}
\newcommand{\bea}{\begin{eqnarray}}
\newcommand{\eea}{\end{eqnarray}}
\newcommand{\Bea}{\begin{eqnarray*}}
\newcommand{\Eea}{\end{eqnarray*}}

\newcommand{\bi}{\begin{itemize}}
\newcommand{\ei}{\end{itemize}}

\newtheorem{Definition}{Definition}[section]
\newtheorem{Theorem}[Definition]{Theorem}
\newtheorem{Lemma}[Definition]{Lemma}

\newtheorem{Corollary}[Definition]{Corollary}
\newtheorem{Conjecture}[Definition]{Conjecture}

\newtheorem{Revised Question}[Definition]{Revised Question}

\theoremstyle{remark}

\begin{document}

\title{A coding theoretic approach to the uniqueness conjecture for projective planes of prime order}
\author{Bhaskar Bagchi\\
Theoretical Statistics and Mathematics Unit\\
Indian Statistical Institute\\
Bangalore - 560 059, India.\\
email: bbagchi@isibang.ac.in}
\date{}

\maketitle



\abstract{An outstanding folklore conjecture asserts that, for any prime $p$, up to isomorphism the projective plane $PG(2,\mathbb{F}_p)$ over the field $\mathbb{F}_p := \mathbb{Z}/p\mathbb{Z}$ is the unique projective plane of order $p$. Let $\pi$ be any projective plane of order $p$. For any partial linear space ${\cal X}$, define the inclusion number $i({\cal X},\pi)$ to be the number of isomorphic copies of ${\cal X}$ in $\pi$. In this paper we prove that if ${\cal X}$ has at most $\log_2 p$ lines, then $i({\cal X},\pi)$ can be written as an explicit rational linear combination (depending only on ${\cal X}$ and $p$) of the coefficients of the complete weight enumerator (c.w.e.) of the $p$-ary code of $\pi$. Thus, the c.w.e. of this code carries an enormous amount of structural information about $\pi$. In consequence, it is shown that if $p > 2^ 9=512$, and $\pi$ has the same c.w.e. as $PG(2,\mathbb{F}_p)$, then $\pi$ must be isomorphic to $PG(2,\mathbb{F}_p)$. Thus, the uniqueness conjecture can be approached via a thorough study of the possible c.w.e. of the codes of putative projective planes of prime order.}

\section{Introduction}

An {\bf incidence system} is a pair $(P,L)$ where $P$ is a set, and $L$ is a set of subsets of $P$. The elements of $P$ and $L$ are called the {\bf points} and {\bf lines} (or blocks) of the incidence system. If $x$ is a point and $\ell$ is a line, then we say that $x$ and $\ell$ are {\bf incident} with each other if $x \in \ell$.

For any two incidence systems ${\cal X}=(P,L)$ and ${\cal X}^\prime =(P^\prime,L^\prime)$, we say that ${\cal X}^\prime$ is a {\bf subsystem} of ${\cal X}$ if $P^\prime \subseteq P$, and every $\ell^\prime \in L^\prime$ can be written as $\ell^\prime=\ell \cap P^\prime$ for some $\ell \in L$.

Let ${\cal X}=(P,L)$ and ${\cal X}^\prime =(P^\prime, L^\prime)$ be two incidence systems. If $f: P \rightarrow P^\prime$ and $g : L \rightarrow L^\prime$ are functions such that, for all $x \in P$ and $\ell \in L$, $x \in \ell \Rightarrow f(x) \in g(\ell)$, then we say that the pair $(f,g)$ is a {\bf homomorphism} from ${\cal X}$ to ${\cal X}^\prime$. If, further, both $f$ and $g$ are one-to-one, then we say that $(f,g)$ is a {\bf monomorphism} from ${\cal X}$ to ${\cal X}^\prime$. If both $f$ and $g$ are onto then $(f,g)$ is called an {\bf epimorphism} from ${\cal X}$ to ${\cal X}^\prime$, and we say that ${\cal X}^\prime$ is a {\bf homomorphic image} of ${\cal X}$. An {\bf isomorphism} from ${\cal X}$ to ${\cal X}^\prime$ is a homomorphism which is both a monomorphism and an epimorphism. If there is an isomorphism from ${\cal X}$ to ${\cal X}^\prime$, then we say that ${\cal X}$ and ${\cal X}^\prime$ are {\bf isomorphic}. One often identifies isomorphic incidence systems. Note that any monomorphism from ${\cal X}$ to ${\cal X}^\prime$ may be viewed as an isomorphism between ${\cal X}$ and a subsystem of ${\cal X}^\prime$. An {\bf automorphism} of ${\cal X}$ is an isomorphism from ${\cal X}$ to ${\cal X}$. The set $\text{Aut}({\cal X})$ of all automorphisms of ${\cal X}$ forms a group (the automorphism group of ${\cal X}$) under componentwise composition.

An incidence system ${\cal X}$ is said to be a {\bf partial linear space} if (a) each line of ${\cal X}$ is incident with at least two points, and (b) any two distinct points of ${\cal X}$ are together incident with at most one line of ${\cal X}$. We say that ${\cal X}$ is a {\bf linear space} if it satisfies (a) and (b$^\prime$): any two distinct points of ${\cal X}$ are together incident with a unique line of ${\cal X}$. Note that if $(f,g)$ is a monomorphism from a partial linear space ${\cal X}$ to a partial linear space ${\cal X}^\prime$, then $g$ is determined by $f$ in view of the requirement $f(\ell) \subseteq g(\ell)~ (\ell \in L)$. Here $f(\ell) := \{f(x) : x\in \ell\}$. In this case, it is usual to omit mention of $g$, and one says that $f$ is the monomorphism. In particular, this remark applies to isomorphisms and automorphisms of partial linear spaces.

If ${\cal X}$ is a partial linear space, then its {\bf dual} ${\cal X}^\ast$ is the incidence system obtained from ${\cal X}$ by interchanging the notions of points and lines. More precisely, the incidence system ${\cal X}^\ast =(P^\ast, L^\ast)$ is the dual of ${\cal X}=(P,L)$ if there are bijections $f: P \rightarrow L^\ast$, $g: L \rightarrow P^\ast$ such that, for all $x\in P$ and $\ell \in L$, $x\in \ell \Leftrightarrow g(\ell) \in f(x)$. This defines the dual up to isomorphisms. Note that the dual ${\cal X}^\ast$ of a partial linear space ${\cal X}$ is again a partial linear space iff each point of ${\cal X}$ is incident with at least two lines of ${\cal X}$. In this case, we have ${\cal X}^{\ast \ast}={\cal X}$.

Finally, a {\bf projective plane} is an incidence system ${\cal X}$ such that (i) ${\cal X}$ is a linear space, (ii) its dual ${\cal X}^\ast$ is also a linear space, and (iii) given any two distinct lines of ${\cal X}$, there is at least one point of ${\cal X}$ which is non-incident with both lines. Note that, in the presence of (i) and (ii), the condition (iii) is equivalent to its dual condition (iii)$^\ast$: given any two distinct points of ${\cal X}$, there is at least one line of ${\cal X}$ which is non-incident with both points. Thus, the dual of a projective plane is a projective plane.

If $x_1,x_2$ are distinct points of a projective plane, we shall denote by $x_1 \vee x_2$ the unique line incident with both $x_1$ and $x_2$. Dually, if $\ell_1, \ell_2$ are distinct lines of a projective plane, then $\ell_1 \wedge \ell_2$ denotes the unique point incident with both $\ell_1$ and $\ell_2$. Note that, for any non-incident point line pair $(x,\ell), ~y \mapsto x \vee y$ defines a function from the set of all points on $\ell$ to the set of all lines through $x$. Clearly the function $m \mapsto m \wedge \ell$ is its inverse. So these two functions are bijections. These bijections are the so-called {\bf perspectivities} in the projective plane. The existence of these perspectivities may be used to see that if $\pi$ is a finite projective plane, then there is a number $n \geq 2$ such that (a) each point of $\pi$ is incident with exactly $n+1$ lines, (b) each line of $\pi$ is incident with exactly $n+1$ points, (c) the total number of points of $\pi$ is $n^2+n+1$, and (d) the total number of lines of $\pi$ is $n^2+n+1$. This number $n$ is called the {\bf order} of the finite projective plane $\pi$.

{\bf Examples} (1) {\bf The field planes}. Let $V$ be a three dimensional vector space over a field $\mathbb{F}$. For $i=1,2$, let $V_i$ be the set of all $i$-dimensional vector subspaces of $V$. We identify each element $\ell$ of $V_2$ with the set of all elements of $V_1$ contained in $\ell$. With this identification, $PG(2,\mathbb{F}) :=(V_1,V_2)$ is a projective plane, called the projective plane over $\mathbb{F}$. With a little care in handling non-commutativity of multiplication, this construction generalizes to yield the projective plane $PG(2,\mathbb{D})$ over any division ring $\mathbb{D}$. Recall that, by a famous theorem of Wedderburn, the finite fields $\mathbb{F}_q$ ($q$ prime power) are the only finite division rings. Specializing the above construction, we get the classical finite projective planes $PG(2,\mathbb{F}_q)$, of order $q$.

(2) {\bf The free projective plane}. The usual definition of the free projective plane may be found in [9]. We like to rephrase this definition slightly, as follows. We first define a sequence ${\cal X}_n =(P_n, L_n), n \geq 1$, of partial linear spaces by induction on $n$. ${\cal X}_1$ is the incidence system whose points and lines are the vertices and edges of the 4-cycle. That is, ${\cal X}_1=(P_1,L_1)$, where $P_1=\{1,2,3,4\}$ and $L_1=\{(1,2\},\{2,3\},\{3,4\},\{4,1\}\}$. Having defined ${\cal X}_n$, we extend it to ${\cal X}_{n+1}$ by introducing a new point $\ell_1 \wedge \ell_2$ corresponding to each pair $\{\ell_1,\ell_2\}$ of lines of ${\cal X}_n$ such that no point of ${\cal X}_n$ is incident in ${\cal X}_{n}$ with both $\ell_1$ and $\ell_2$, and introducing a new line $x_1 \vee x_2$ corresponding to each pair $\{x_1,x_2\}$ of points of ${\cal X}_n$ such that no line of ${\cal X}_n$ is incident in ${\cal X}_n$ with both $x_1$ and $x_2$. The point $\ell_1 \wedge \ell_2$ is incident in ${\cal X}_{n+1}$ with the lines $\ell_1, \ell_2$ and with no other line. The line $x_1 \vee x_2$ is incident in ${\cal X}_{n+1}$ with the points $x_1,x_2$ and with no other point. The incidences between the old points and old lines are as in ${\cal X}_n$. Thus, by construction, each ${\cal X}_n$ is a subsystem of ${\cal X}_{n+1}$. Finally, we define ${\cal F} =\left(\bigcup\limits^\infty_{n=1} P_n, \bigcup\limits^\infty_{n=1} L_n\right)$. Clearly ${\cal F}$ is an (infinite) projective plane. It is called the free projective plane. An easy induction shows that each ${\cal X}_n$ is {\bf self-dual} (i.e., isomorphic to its dual). It follows that (like the field planes) the free projective plane is self dual.

A projective plane $\sigma$ is said to be a {\bf subplane} of the projective plane $\pi$ if $\sigma$ is a subsystem of $\pi$. A projective plane $\pi$ is said to be {\bf prime} if it has no proper projective subplane. The projective planes over prime fields $\mathbb{F}$ (i.e., $\mathbb{F}=\mathbb{Q}$ or $\mathbb{F}=\mathbb{F}_p$ for some prime $p$) are obvious examples of prime projective planes. It is not hard to see that the free projective plane ${\cal F}$ is also prime. (Indeed, using the above construction, one sees that the automorphism group of ${\cal F}$ is transitive on the 4-cycles - isomorphic copies of ${\cal X}_1$ - in ${\cal F}$. The same is true of $PG(2,\mathbb{F}), \mathbb{F}$ a prime field). In a private conversation with the author some years back, N.M. Singhi forwarded:

\begin{Conjecture} (Singhi). The projective planes over prime fields and the free projective plane are the only examples of prime projective planes.
\end{Conjecture}

One of the reasons why the free projective plane is important is the fact that every prime projective plane is a homomorphic image of ${\cal F}$. (Indeed, given a prime projective plane $\pi$, one may readily find a monomorphism $f_1$ of ${\cal X}_1$ into $\pi$. Hence one may inductively construct a homomorphism $f_n$ of ${\cal X}_n$ into $\pi$ such that $f_{n+1}$ extends $f_n$ for each $n$. Then $f:= \bigcup\limits_{n \geq 1} f_n$ is a homomorphism of ${\cal F}$ into $\pi$. Since $\pi$ is prime, $f$ must be an epimorphism.) This indicates that an in-depth study of the partial linear spaces ${\cal X}_n$ may be an useful approach to Singhi's conjecture. Note that, in particular, if $\pi$ is a finite prime projective plane, then it follows from the above that $\pi$ is a homomorphic image of ${\cal X}_n$ for sufficiently large $n$. We suggest that an investigation of the codes of ${\cal X}_n$ (over various primes) may be fruitful. But we do not even know a formula for the number $v_n$ of points (= number of lines) of ${\cal X}_n$. The sequence $\{v_n\}$ begins with $4,6,7,9,13,33,\ldots$.

When $q$ is a ``genuine" prime power (i.e., $q=p^e, p$ prime, $e \geq 2$) and $q>8$, there are many constructions of projective planes of order $q$ which are not field planes. But no such construction is known for $q=p$. The subject of this paper is the following

\begin{Conjecture}\label{1.2} (folklore): Up to isomorphism, $PG(2,\mathbb{F}_p)$ is the only projective plane of prime order $p$.
\end{Conjecture}

Conjecture \ref{1.2} bears a superficial resemblance to Singhi's conjecture. But the relationship between these two is far from clear. We do not even know if a projective plane of prime order must ba a prime projective plane, or if a prime projective plane which is finite must be of prime order. Another related conjecture is due to H. Neumann \cite{11} (which is much stronger than the finite case of Singhi's conjecture): she conjectured that a finite projective plane has no projective subplane of order two (if and) only if it is isomorphic to $PG(2,\mathbb{F}_q)$ for some odd prime power $q$.

In the humble opinion of this author, the uniqueness conjecture \ref{1.2} is one of the most beautiful and important open problems in mathematics. It is amusing as well as sad that it finds no mention in the lists of ``problems for the new millenium" compiled by various authors at the turn of the century. It is not lacking in history and pedigree. With some imagination, one may trace the history of such problems in finite geometry back to Euler's 1782 paper \cite{7} on the problem of the thirty six officers. The projective planes $PG(2,\mathbb{F}_p)$ were first constructed by von Staudt \cite{12} in 1856. They were generalized to the planes $PG(2,\mathbb{F}_q),q$ prime power, by Fano \cite{9} in 1892. The first examples of non-field finite projective planes were constructed by Veblen and Wedderburn \cite{13} in 1907. The conjecture \ref{1.2} must have occurred to these early authors. The vast literature on finite projective planes include (usually as special cases) many characterizations of $PG(2,\mathbb{F}_p)$ on the assumption of moderately large automorphism groups. See \cite{10} and
\cite[Chapter~5]{6} for many of these results.

\section{Coding theory}

One very fruitful approach to problems in finite geometry has been through the study of codes attached to these geometries. For a comprehensive account of these connections, the reader may consult \cite{1}.

If $P$ is a finite set and $p$ is a prime number, then consider the $\mathbb{F}_p$-vector space $\mathbb{F}_p^P$ consisting of all functions $f: P \rightarrow \mathbb{F}_p$. For any such $f$, the {\bf support} of $f$ is the set $\{x \in P : f(x) \neq 0\}$, and the {\bf Hamming weight} $|f|$ of $f$ is the size (cardinality) of the support of $f$. The {\bf type} of $f$ is the $p$-tuple $(a_\alpha : \alpha \in \mathbb{F}_p$) where $a_\alpha =\# \{x \in P : f(x)=\alpha\}$. $\mathbb{F}_p^P$ is a metric space with the Hamming metric given by $d(f,g)=|f-g|$. One also equips $\mathbb{F}_p^P$ with the usual non-degenerate symmetric bilinear form $\langle \cdot, \cdot \rangle$ defined by $\langle f,g\rangle =\sum\limits_{x\in P} f(x)g(x)$.

A $p$-ary code $C$ is a linear subspace of $\mathbb{F}_p^P$, viewed as a metric space with the Hamming metric inherited from $\mathbb{F}_p^P$. The vectors in $C$ are called the {\bf words} of $C$. The elements of $P$ are called the {\bf co-ordinate positions} of $C$. The {\bf minimum weight} of $C$ is defined to be the number $\min\{|f| : f\in C \backslash \{0\}\}$. The {\bf dual code} $C^\perp$ of $C$ is the orthocomplement of $C$ with respect to $\langle \cdot, \cdot \rangle$. Thus $C^\perp := \{g \in \mathbb{F}_p^P : \langle g,f\rangle =0 ~\forall ~f \in C\}$. The {\bf Hamming weight enumerator} of $C$ is the polynomial $F\in \mathbb{Z}[X,Y]$ given by $F(X,Y)=\sum\limits_{f \in C} X^{n-|f|} Y^{|f|}$, where $n:= \# (P)$ is the so-called length of $C$. If $\underline{Z}=(Z_\alpha : \alpha \in \mathbb{F}_p)$ is a $p$-tuple of commuting variables, then the {\bf complete weight enumerator} $G(\underline{Z})$ of $C$ is defined to be the polynomial $G(\underline{Z})=\sum\limits_{f\in C} \underline{Z}^{\text{type}(f)}$. Here we have used the usual notation $\underline{Z}^{\underline{i}}$ for $\Pi_\alpha Z_\alpha^{i_\alpha}$, where $\underline{i} =(i_\alpha : \alpha \in \mathbb{F}_p)$ is a multi-index. Note that the complete weight enumerator (c.w.e.) of $C$ carries much more information about the code than the Hamming weight enumerator. Indeed, $F$ may be obtained from $G$ by the substitutions $X=Z_0, Y=Z_1=Z_2=\ldots=Z_{p-1}$. The Hamming weight enumerator enumerates the frequencies of various Hamming weights occurring in the code, while the complete weight enumerator enumerates the frequencies of the various types. Finally, note that - since $\langle \cdot, \cdot \rangle$ is non-degenerate - we have the usual formula relating the dimensions of $C$ and $C^\perp : \dim(C) +\dim (C^\perp)=n$. There are also beautiful formulae giving the Hamming weight enumerator (and, more generally, the c.w.e.) of $C^\perp$ in terms of the corresponding enumerator of $C$. While these formulae are extremely important, we shall have no occasion to use them in this paper.

Let ${\cal X}$ be a finite incidence system and $p$ be a prime. Let $P$ be the point set of ${\cal X}$. For any line $\ell$ of ${\cal X}$, we consider its {\bf indicator function} $\ell : P \rightarrow \mathbb{F}_p$ given by $\ell(x) =1$ if $x\in \ell$, and $\ell(x) =0$ if $x\not\in \ell$. Note that we have used the same letter to denote a line and its indicator. Thus, lines of ${\cal X}$ are also code words in $\mathbb{F}^P_p$. {\bf The $p$-ary code $C_p({\cal X})$ of ${\cal X}$} is defined to be the vector subspace of $\mathbb{F}_p^P$ spanned by all the lines of ${\cal X}$.

If $\pi$ is a finite projective plane of order $n$, then it is easy to see that $C_p(\pi)$ is trivial when $p$ does not divide $n$. However, when $p$ divides $n$, $C_p(\pi) \cap C^\perp_p(\pi)$ is of co-dimension one in $C_p(\pi)$ (this intersection, the so-called ``hull", is actually spanned by the pairwise differences of the lines of $\pi$). Also, when $p$ divides $n$, the minimum weight of $C_p(\pi)$ is $n+1$, and the minimum weight words of $C_p(\pi)$ are precisely the non-zero scalar multiples of the lines of $\pi$. When $p$ exactly divides $n$ (i.e., $p \mid n$ but $p^2 \nmid n$) $C^\perp_p(\pi)$ is a subcode (of codimension one) in $C_p(\pi)$, so that $\dim (C^\perp_p (\pi))=\binom{n+1}{2}$ and $\dim (C_p(\pi))=\binom{n+1}{2}+1$. See \cite{1} for the proofs of these results.

In particular, if $\pi$ is a projective plane of prime order $p$, $C_p^\perp (\pi) \subset C_p(\pi)$, $\dim C^\perp_p (\pi)=\binom{p+1}{2}$, $\dim C_p(\pi)=\binom{p+1}{2}+1$, and the minimum weight words of $C_p(\pi)$ are the non-zero scalar multiples of lines. Also, Inamdar \cite{10a} proved that, in this case, the minimum weight words of $C^\perp_p(\pi)$ are precisely the non-zero scalar multiples of the pairwise differences of lines of $\pi$.

These results apply, in particular, to the $p$-ary code of $PG(2,\mathbb{F}_p)$. In \cite{4}, the present author proved that the first four minimum weights of the $p$-ary code of $PG(2,\mathbb{F}_p) ~ (p \geq 5)$ are $p+1,2p,2p+1,$ and $3p-3$. In an earlier paper \cite{8}, Fack et al had proved that the only words of weight $p+1,2p$ or $2p+1$ (in this particular code) are the non-zero $\mathbb{F}_p$-linear combinations of pairs of lines in $PG(2,\mathbb{F}_p)$. However, a complete classification of the words of Hamming weight $3p-3$ in $C_p(PG(2,\mathbb{F}_p))$ remains an open problem. We shall return to this question in the concluding section of this paper.

\section{A series of Lemmas}

In this section we prove a number of lemmas culminating in Lemma \ref{3.9} which will be used in the next section. Our first lemma is well known.

\begin{Lemma} \label{3.1} Let $P$ be the point-set of a projective plane $\pi$ of prime order $p$. Let $w \in \mathbb{F}^P_p$. Then $w \in C_p(\pi)$ if and only if $\langle w,\ell \rangle =\langle w,\mathbf{1}\rangle$ for all lines $\ell$ of $\pi$. (Here $\mathbf{1}$ is the constant function 1 on $P$.)
\end{Lemma}

{\bf Proof:} This means that $w \in C_p$ iff $w$ is orthogonal to all the words in the set $\{\mathbf{1}-\ell:\ell$ a line of $\pi\}$. This is true since this set spans $C^\perp_p$.
\hfill${\Box}$

\begin{Lemma} \label{3.2} Let ${\cal X}$ be a finite partial linear space and let $p$ be a prime. Then $\dim (C_p ({\cal X}^\ast))=\dim (C_p({\cal X}))$.
\end{Lemma}


{\bf Proof:} Let ${\cal X}=(P,L),~ v =\# (P), ~b=\#(L)$. Let $N$ be the $v\times b$ incidence matrix of ${\cal X}$. That is, the rows and columns of $N$ are indexed by $P$ and $L$ respectively, and - for $x\in P, \ell \in L$ - the $(x,\ell)$th entry of $N$ is $=1$ if $x\in \ell$ and $=0$ otherwise. If we view $N$ as a linear operator from $\mathbb{F}_p^L$ to $\mathbb{F}_p^P$, then $C_p({\cal X})$ is precisely the image of $N$. Therefore $\dim C_p({\cal X})=\text{rank}_p(N)$. Now note that the transposed matrix $N^\ast$ is the incidence matrix of ${\cal X}^\ast$, so that we also have $\dim C_p({\cal X}^\ast)=\text{rank}_p(N^\ast)$. As $\text{rank}_p(N^\ast)=\text{rank}_p(N)$, the result follows. \hfill${\Box}$

\begin{Lemma}\label{3.3} Let $\pi$ and $\sigma$ be two projective planes of prime order $p$. Suppose $\pi$ and $\sigma$ share at least $p^2+1$ lines. Then $\pi=\sigma$.
\end{Lemma}

{\bf Proof:} Let $L_0$ be a set of $p^2+1$ lines common to $\pi$ and $\sigma$. Since $\pi$ has $p^2+p+1$ lines, and each point of $\pi$ is in $p+1$ lines, it follows that the union of the lines in $L_0$ is the entire point set of $\pi$. Similarly for $\sigma$. So $\pi$ and $\sigma$ have the same point set, call it $P$. Let $C_0$ be the subcode of $\mathbb{F}_p^P$ spanned by $L_0$. Thus $C_0$ is a subcode of both $C_p(\pi)$ and $C_p(\sigma)$.

Consider the incidence system ${\cal X}=(P,L_0)$. Thus, $C_0=C_p({\cal X})$. Consider the restriction map $\rho : \mathbb{F}_p^L \rightarrow \mathbb{F}_p^{L_0}$ given by $w \mapsto w \mid_{L_0}$, where $L$ is the full set of lines of $\pi$. $\rho$ is a linear map which restricts to a linear map from $C_p(\pi^\ast)$ onto $C_p({\cal X}^\ast)$. The kernel of this restricted map consists of the words $w$ of $C_p(\pi^\ast)$ with support$(w) \subseteq L \backslash L_0$. But $\#(L \backslash L_0)=p$ and $C_p(\pi^\ast)$ has no non-zero word of Hamming weight $\leq p$. Therefore, the kernel is trivial and so $\rho$ restricts to a vector space isomorphism between $C_p(\pi^\ast)$ and $C_p({\cal X}^\ast)$. In conjunction with Lemma \ref{3.2}, this yields $\dim(C_p(\pi))=\dim(C_p(\pi^\ast))=\dim(C_p({\cal X}^\ast))=\dim(C_p({\cal X}))$. Since $C_0=C_p({\cal X})$ is a subcode of $C_p(\pi)$, it follows that $C_p(\pi)=C_0$. Similarly, $C_p(\sigma)=C_0$. Thus, $C_p(\pi)=C_p(\sigma)$. But the lines of $\pi$ are precisely the supports of the minimum weight words of $C_p(\pi)$, and similarly for $\sigma$. Therefore $\pi$ and $\sigma$ have the same set of lines as well. Hence $\pi=\sigma$. \hfill${\Box}$

\begin{Definition} \label{3.4} We shall say that an incidence system ${\cal Y}$ is $\mathbf{p}$-{\bf admissible} if it satisfies (i) ${\cal Y}$ has exactly $p^2+p+1$ points, (ii) each line of ${\cal Y}$ is incident with exactly $p+1$ points, and (iii) any two distinct lines of ${\cal Y}$ are together incident with a unique point.
\end{Definition}

Thus, any $p$-admissible incidence system is a partial linear space. Note that any set of lines in a projective plane of order $p$ may be viewed as the set of all lines of a $p$-admissible incidence system.

\begin{Lemma} \label{3.5} Let $p \geq 2$. Let $\sigma$ be a $p$-admissible system. Then $\sigma$ has at most $p^2+p+1$ lines. Equality holds iff $\sigma$ is a projective plane of order $p$.
\end{Lemma}

{\bf Proof:} Fix a point $x$ of $\sigma$. Then the lines of $\sigma$ through $x$, minus the point $x$, are pairwise disjoint subsets of size $p$ each in the set of $p(p+1)$ remaining points. So each of the $p^2+p+1$ points $x$ is in at most $p+1$ lines. But each of the lines of $\sigma$ contains exactly $p+1$ points. So, a two-way counting shows that $\sigma$ has at most $p^2+p+1$ lines.

Now suppose $\sigma$ has $p^2+p+1$ lines. Then the above argument shows that each point $x$ is in exactly $p+1$ lines. Therefore the lines through $x$ induce a partition of the remaining points. So, if $y \neq x$ is another point, then a unique line joins $x$ and $y$. Since $x$ was arbitrary, this shows that any two distinct points of $\sigma$ are together in a unique line of $\sigma$.

Now fix two lines $\ell_1 \neq \ell_2$ of $\sigma$. Let $x$ be the common point of $\ell_1$ and $\ell_2$. Since $p+1 \geq 3$, there is a third line $\ell$ through $x$ and there is a point $y \neq x$ on $\ell$. Then $y$ is non-incident with both $\ell_1$ and $\ell_2$. So $\sigma$ is a projective plane of order $p$. \hfill${\Box}$

\begin{Lemma} \label{3.6} Let $S$ be the union of $k \geq 1$ lines of a $p$-admissible incidence system. Then $(p+1)k-\binom{k}{2} \leq \# (S) \leq pk+1$.
\end{Lemma}

{\bf Proof:} Induction on $k$. The result is trivial for $k=1$. So let $k>1$, and let $S=\bigcup\limits_{0 \leq i < k} \ell_i$, where $\ell_i$'s are distinct lines of the $p$-admissible incidence system. Write $S=S^\prime \cup \ell_{k-1}$, where $S^\prime =\bigcup\limits_{0 \leq i < k-1} \ell_i$.

Since any two distinct lines meet at a unique point, we have $1 \leq \# (S^\prime \cap \ell_{k-1}) \leq k-1$. Since, also, each line is of size $p+1$, we get
\begin{eqnarray*}
k(p+1)-\binom{k}{2} &=& (k-1)(p+1) -\binom{k-1}{2} +(p+1)-(k-1) \\
&\leq& \# (S^\prime) + \#(\ell_{k-1}) -\#(S^\prime \cap \ell_{k-1})\\
&\leq& (k-1)p+1+(p+1)-1\\
&=& kp+1.
\end{eqnarray*}
As $\#(S^\prime)+\#(\ell_{k-1})-\# (S^\prime \cap \ell_{k-1})=\#(S)$, this completes induction. \hfill${\Box}$

\begin{Lemma} \label{3.7} Let ${\cal Y}$ and ${\cal Y}^\prime$ be two $p$-admissible incidence systems. Suppose the union of some $m$ lines of ${\cal Y}$ equals the union of some $k$ lines of ${\cal Y}^\prime$. If $\binom{k}{2} < p$ then $m=k$.
\end{Lemma}

{\bf Proof:} Let $S$ be a set which is the union of $m$ lines of ${\cal Y}$ as well as the union of $k$ lines of ${\cal Y}^\prime$. Since $p > \binom{k}{2}$, Lemma \ref{3.6} implies that $(k-1)p+1< (p+1)k-\binom{k}{2} \leq \# (S) \leq mp +1$. Therefore $m \geq k$. Suppose, if possible, that $m >k$. Then $S$ is the union of $k$ lines of ${\cal Y}^\prime$ and $S$ contains the union of $k+1$ lines of ${\cal Y}$. Therefore, Lemma \ref{3.6} implies that we have
$(k+1)(p+1)-\binom{k+1}{2} \leq \# (S) \leq kp+1$. Hence $p \leq \binom{k}{2}$. But this contradicts our assumption. \hfill${\Box}$

\begin{Lemma}\label{3.8} Let $k$ be a positive integer, and let $x_i, 0 \leq i < k$, be $k$ non-negative integers such that $\sum\limits_{0 \leq i < k} 2^i x_i =2^k-1$. Then,
\begin{enumerate}
\item[(a)] $\sum\limits_{0 \leq i <k} x_i \geq k$, and
\item[(b)] if $\sum\limits_{0 \leq i < k} x_i=k$ then $x_i=1$ for all $i$.
\end{enumerate}
\end{Lemma}

{\bf Proof:} It suffices to show that if $\sum_{0 \leq i  < k} x_i \leq k$ then $x_i=1$ for all $i$. We prove this by induction on $k$. The result is trivial for $k=1$. So assume $k >1$. Note that $\sum\limits_{0 \leq i < k} 2^i x_i =2^k-1$ implies that $x_0$ is odd. In particular, $x_0 \geq 1$. We define $k-1$ non-negative integers $y_i, 0 \leq i < k-1$, as follows. $y_0 =\frac{1}{2} (x_0-1)+x_1$, and $y_i = x_{i+1}$ for $1 \leq i < k-1$. We have $\sum\limits_{0 \leq i < k-1} 2^i y_i =\frac{1}{2} (-1+\sum\limits_{0 \leq i < k} 2^i x_i)=2^{k-1}-1$, and $\sum\limits_{0 \leq i < k-1} y_i< \sum\limits_{0 \leq i < k} x_i \leq k$. Therefore, the induction hypothesis implies that $y_i=1$ for $0\leq i < k-1$. That is, $x_i =1$ for $ 1<i < k$, and $x_0+2x_1=3$. So, either $x_0=x_1=1$, or $x_0=3,x_1=0$. But, in the latter case, we get $\sum\limits_{0 \leq i < k} x_i=k+1$, contrary to our assumption. So $x_i=1$ for all $i$. This completes induction. \hfill${\Box}$

\begin{Lemma} \label{3.9} Let $p$ be a prime, and ${\cal Y}$ be a $p$-admissible incidence system with exactly $k$ lines. Enumerate the lines of ${\cal Y}$ (in any order) as $\ell_i, ~0 \leq i < k$. Consider the word $w \in C_p({\cal Y})$ given by $w=\sum\limits_{0 \leq i < k}2^i \ell_i$. Let $\pi$ be a projective plane of order $p$, and suppose $w^\prime$ is a word in $C_p(\pi)$ such that $\text{type} (w^\prime)=\text{type}(w)$. If $p \geq 2^k$, then there are lines $\ell_i^\prime, ~0 \leq i < k$, of $\pi$ such that $w^\prime =\sum\limits_{0 \leq i < k} 2^i \ell^\prime_i$. Further, there is a monomorphism $f$ from ${\cal Y}$ into $\pi$ such that $f(\ell_i) =\ell_i^\prime, ~0 \leq i < k$.
\end{Lemma}

{\bf Proof:} For integers $i \geq 0$ and $x\geq 0$, let $\delta_i(x)$ denote the $i^{\text{th}}$ digit in the binary expansion of $x$, counting from the right, and taking the rightmost digit as the $0^{\text{th}}$.

Let $P$ and $Q$ be the point sets of $\pi$ and ${\cal Y}$, respectively. For $0 \leq i <k$, define
\[
\ell^\prime_i =\{x \in P : \delta_i (w^\prime (x))=1\}.
\]
(Here we have identified the elements of $\mathbb{F}_p$ with the integers $0,1,\ldots, p-1$.) Notice that, since $w=\sum\limits_{0 \leq i < k} 2^i \ell_i$, and $p \geq 2^k$, we also have
\[
\ell_i =\{ x \in Q : \delta_i (w(x))=1\}
\]
for $ 0 \leq i < k$. Further, by the definition of $\ell_i^\prime$, we have
\[
w^\prime =\sum\limits_{0 \leq i < k} 2^i \ell_i^\prime.
\]
As $\text{type}(w^\prime) =\text{type}(w)$, there is a bijection $f: Q \rightarrow P$ such that $w=w^\prime \circ f$. It follows that for $x\in Q$, and $0 \leq i < k$,
\[
x \in \ell_i \Leftrightarrow \delta_i (w(x))=1 \Leftrightarrow \delta_i (w^\prime (f(x)))=1 \Leftrightarrow f(x) \in \ell^\prime_i.
\]
Thus $f(\ell_i)=\ell^\prime_i$ for $0 \leq i < k$. Therefore, if ${\cal Y}^\prime$ denotes the incidence system with point set $P$ and lines $\ell_i^\prime, ~0 \leq i < k$, then $f$ is an isomorphism between ${\cal Y}$ and ${\cal Y}^\prime$. In consequence, ${\cal Y}^\prime$ is also $p$-admissible, so that $\ell_i^\prime$ are sets of size $p+1$ each, and any two distinct $\ell_i^\prime$'s meet at a unique point.

Thus, to complete the proof, we have to show that the $k$ sets $\ell_i^\prime$ are lines of $\pi$. So far we have not used the assumption $w^\prime \in C_p(\pi)$. This assumption will play a crucial role in what follows.

Let $S$ and $S^\prime$ be the supports of $w$ and $w^\prime$ respectively. Since $p \geq 2^k$, we have $S=\bigcup\limits_{0 \leq i < k} \ell_i$ and $S^\prime =\bigcup\limits_{0 \leq i < k} \ell_i^\prime$.

{\bf Claim 1:} If $\ell$ is a line of $\pi$ such that $\ell \nsubseteq S^\prime$, then $\sum\limits_{x \in \ell} w^\prime (x)=2^k-1$ (equality in $\mathbb{N}$).

To prove this claim, first we note that
\begin{eqnarray*}
\sum\limits_{x \in P}w^\prime (x) &\equiv& \langle w^\prime, \mathbf{1}\rangle \equiv \langle w, \mathbf{1}\rangle \equiv \left\langle \sum\limits_{0 \leq i < k} 2^i \ell_i, \mathbf{1} \right\rangle \\
&\equiv& \sum\limits_{0 \leq i < k} 2^i \equiv 2^k-1 (\text{mod}~p).
\end{eqnarray*}
Since $w^\prime \in C_p(\pi)$, Lemma \ref{3.1} implies that $\sum\limits_{x \in \ell} w^\prime (x)\equiv 2^k -1 (\text{mod}~p)$ for any line $\ell$ of $\pi$. Therefore we have, as $p \geq 2^k$,
\begin{equation}
\sum\limits_{x \in \ell} w^\prime (x) \geq 2^k -1
\end{equation}
for any line $\ell$ of $\pi$ (inequality in $\mathbb{N}$).

Now fix a point $y \not\in S^\prime$. Adding the inequalities (1) over all the $p+1$ lines $\ell$ through $y$, we get
\begin{eqnarray*}
(p+1) (2^k-1) &\leq& \sum\limits_{\ell \ni y} \sum\limits_{x\in \ell} w^\prime (x) \\
&=& \sum\limits_{x\in P} w^\prime (x) \\
&=& \sum\limits_{x\in P} \sum\limits_{0 \leq i < k} 2^i \ell_i^\prime (x)\\
&=& \sum\limits_{0 \leq i < k} 2^i \sum\limits_{x\in P} \ell_i^\prime (x) \\
&=& (p+1)(2^k-1).
\end{eqnarray*}
Since the two extreme terms here are equal, we must have equality throughout this argument. Therefore we have equality in (1) for any line $\ell$ through $y$. Since $y \not\in S^\prime$ was an arbitrary point, we have equality for any line $\ell \nsubseteq S^\prime$. This proves Claim 1.

{\bf Claim 2:} For any line $\ell \nsubseteq S^\prime$ of $\pi$, we have $\# (\ell \cap \ell^\prime_i)=1$ for $0 \leq i < k$.

To prove this claim, note that
\begin{eqnarray*}
\sum\limits_{0 \leq i < k} 2^i \# (\ell \cap \ell^\prime_i) &=& \sum\limits_{x \in \ell} \sum\limits_{0 \leq i < k} 2^i \ell^\prime_i (x) \\
&=& \sum\limits_{x \in \ell} w^\prime (x)\\
&=& 2^k-1 ~(\text{by ~Claim~1}).
\end{eqnarray*}
Therefore, Lemma \ref{3.8} (a) implies that
\begin{equation}
\sum\limits_{0 \leq i < k} \# (\ell \cap \ell_i^\prime) \geq k
\end{equation}
for any line $\ell \nsubseteq S^\prime$.

Again, fix a point $y \not\in S^\prime$, and add the inequality (2) over all $p+1$ lines $\ell$ of $\pi$ through $y$. We get:
\begin{eqnarray*}
(p+1) k &\leq& \sum\limits_{\ell \ni y} \sum\limits_{0 \leq i < k} \# (\ell \cap \ell_i^\prime)\\
&=& \sum\limits_{0 \leq i < k} \sum\limits_{\ell \ni y} \# (\ell \cap \ell^\prime_i)\\
&=& \sum\limits_{0 \leq i < k} \# (\ell_i^\prime)\\
&=& (p+1)k.
\end{eqnarray*}
Since the two extreme terms here are equal, we must have equality throughout this argument. Thus, we have equality in (2) for any line $\ell$ through $y$. Since the point $y \not\in S^\prime$ was arbitrary, it follows that we have equality in (2)  for any line $\ell \nsubseteq S^\prime$. Therefore Lemma \ref{3.8} (b) implies that $\# (\ell \cap \ell^\prime_i)=1$ for $0 \leq i < k$ and for any such line $\ell$. This proves Claim 2.

{\bf Claim 3:} $S^\prime$ contains exactly $k$ lines of $\pi$.

To see this, let $m$ be the number of lines $\ell$ of $\pi$ such that $\ell \subseteq S^\prime$. Since $S^\prime$ is the union of the sets $\ell_i^\prime ~(0 \leq i < k)$, and since, by Claim 2, for any two points $x \neq y$ in $\ell_i^\prime$ the line $\ell$ of $\pi$ joining $x$ and $y$ is contained in $S^\prime$, it follows that $S^\prime$ is the union of $m$ lines of $\pi$, as well as the union of $k$ lines of ${\cal Y}^\prime$. Since $p \geq 2^k > \binom{k}{2}$, Lemma \ref{3.7} implies that $m=k$. This proves Claim 3.

Now, let $\sigma$ be the incidence system obtained from $\pi$ by deleting the $k$ lines contained in $S^\prime$ and replacing them by the $k$ lines of ${\cal Y}^\prime$ (contained in $S^\prime$). Since $\pi$ is a projective plane of order $p$ and ${\cal Y}^\prime$ is $p$-admissible, Claim 2 implies that $\sigma$ is $p$-admissible. Since $\sigma$ has $p^2+p+1$ lines, Lemma \ref{3.5} implies that $\sigma$ is also a projective plane of order $p$. Also, by construction, $\sigma$ and $\pi$ share at least $p^2+p+1-k \geq p^2+p+1 - \log_2 p \geq p^2+1$ lines. Therefore, by Lemma \ref{3.3}, $\sigma=\pi$. Since $\ell^\prime_i ~(0 \leq i < k)$ are lines of $\sigma$, it follows that they are lines of $\pi$. \hfill${\Box}$

\section{The main results}

We now introduce:

{\bf Notation:} Let ${\cal Y}$ and ${\cal X}$ be any two finite incidence systems. Then $I({\cal Y},{\cal X})$ will denote the number of monomorphisms from ${\cal Y}$ into ${\cal X}$. Also, $i({\cal Y},{\cal X})$ will denote the number of isomorphic copies of ${\cal Y}$ which are subsystems of ${\cal X}$.

\begin{Lemma} \label{4.1} For any two finite incidence systems ${\cal Y}$ and ${\cal X}$, we have $I({\cal Y},{\cal X})=\# (\text{Aut} ({\cal Y})) \cdot i ({\cal Y},{\cal X})$.
\end{Lemma}

{\bf Proof:} Note that, for any monomorphism $f$ from ${\cal Y}$ to ${\cal X}$, the image ${\cal Y}^\prime$ of ${\cal Y}$ under $f$ is an isomorphic copy of ${\cal Y}$ in ${\cal X}$, and $f$ may be viewed as an isomorphism from ${\cal Y}$ to ${\cal Y}^\prime$. Conversely, for any isomorphic copy ${\cal Y}^\prime$ of ${\cal Y}$ in ${\cal X}$, any isomorphism from ${\cal Y}$ to ${\cal Y}^\prime$ may be viewed as a monomorphism from ${\cal Y}$ to ${\cal X}$. Therefore, to complete the proof, it suffices to show that, whenever ${\cal Y}$ and ${\cal Y}^\prime$ are isomorphic finite incidence systems, the number of isomorphisms from ${\cal Y}$ to ${\cal Y}^\prime$ equals $\#(\text{Aut}({\cal Y}))$. To see this, fix any isomorphism $f$ from ${\cal Y}$ to ${\cal Y}^\prime$, and note that $g \mapsto f \circ g$ is a bijection from the set of all automorphisms of ${\cal Y}$ onto the set of all isomorphisms from ${\cal Y}$ to ${\cal Y}^\prime$. \hfill${\Box}$

{\bf Notation:} For any prime $p$, let ${\cal J}_p$ denote the set of all multi-indices $\underline{j}=(j_\alpha : \alpha \in \mathbb{F}_p)$ such that $|\underline{j}|:= \sum\limits_{\alpha \in \mathbb{F}_p} j_\alpha =p^2+p+1$. Also, let $\underline{X}=(X_\alpha : \alpha \in \mathbb{F}_p)$ be a set of commuting variables.

\begin{Theorem} \label{4.2} Let $\pi$ be a projective plane of prime order $p$, and let $f(\underline{X})=\sum\limits_{\underline{j} \in {\cal J}_p} a_{\underline{j}} \underline{X}^{\underline{j}}$ be the complete weight enumerator of $C_p(\pi)$. (Thus, for $\underline{j} \in {\cal J}_p, ~a_{\underline{j}}$ is the number of words of type $\underline{j}$ in $C_p(\pi)$.) Then, for any partial linear space ${\cal X}$ with at most $\log_2 p$ lines, there are rational numbers $\alpha_{\underline{j}}, ~\underline{j} \in {\cal J}_p$, depending only on ${\cal X}$ and $p$, such that $i({\cal X},\pi)=\sum\limits_{\underline{j}\in {\cal J}_p} \alpha_{\underline{j}} a_{\underline{j}}$.
\end{Theorem}

{\bf Proof:} Let $k$ be the number of lines of ${\cal X}$. Thus $p \geq 2^k$. Notice that, up to isomorphism, there are only finitely many $p$-admissible incidence systems ${\cal Y}$, with exactly $k$ lines, such that ${\cal X}$ is a subsystem of ${\cal Y}$. Let ${\cal Y}_j, ~0 \leq j < m$, be mutually non-isomorphic incidence systems such that every such incidence system ${\cal Y}$ is isomorphic to exactly one ${\cal Y}_j$.

Note that, for any isomorphic copy ${\cal X}^\prime$ of ${\cal X}$ in $\pi$, there is a unique index $j$, $0 \leq j < m$, and a unique isomorphic copy $Y^\prime_j$ in $\pi$ of the incidence system ${\cal Y}_j$, such that ${\cal X}^\prime$ is a subsystem of ${\cal Y}^\prime_j$. (Since ${\cal X}^\prime$ is a partial linear space which is a subsystem of $\pi$, each line $\ell$ of ${\cal X}^\prime$ is contained in a unique line $\overline{\ell}$ of $\pi$. Then ${\cal Y}^\prime_j$ must be the unique subsystem of $\pi$ such that the point set of ${\cal Y}^\prime$ equals that of $\pi$, and the lines of ${\cal Y}^\prime_j$ are the lines $\overline{\ell}$ of $\pi$ as $\ell$ varies over the $k$ lines of ${\cal X}^\prime$.) Therefore we have $i({\cal X},\pi)=\sum\limits_{0 \leq j < m} i({\cal X},{\cal Y}_j)i({\cal Y}_j,\pi)$. Hence, to complete the proof, it suffices to show that for each index $j$, $i({\cal Y}_j,\pi)$ can be written as a rational linear combination of the coefficients of $f$, with coefficients depending only on ${\cal Y}_j$.

So, fix a $p$-admissible incidence system ${\cal Y}$ with $k$ lines, $p \geq 2^k$. We have to show that there are rational numbers $\beta_{\underline{i}}, ~\underline{i} \in {\cal J}_p$, depending only on ${\cal Y}$, such that $i({\cal Y},\pi)=\sum\limits_{\underline{i} \in {\cal J}_p} \beta_{\underline{i}} a_{\underline{i}}$. To see this, take a word $w \in C_p({\cal Y})$ as defined in Lemma \ref{3.9}. Let $\underline{j} \in {\cal J}_p$ be the type of $w$. Then, for each monomorphism $f$ from ${\cal Y}$ to $\pi$, $w \circ f^{-1}$ is one of the $a_{\underline{j}}$ words of type $\underline{j}$ in $C_p(\pi)$. Conversely, if $w^\prime \in C_p(\pi)$ is a word of type $\underline{j}$, then, by the proof of Lemma \ref{3.9}, each of the $\underline{j}!$ bijections $f$ (from the point set of ${\cal Y}$ to the point set of $\pi$) satisfying $w^\prime=w \circ f^{-1}$ is a monomorphism from ${\cal Y}$ to $\pi$. Thus $I({\cal Y},\pi)=\underline{j} ! a_{\underline{j}}$, and hence by Lemma \ref{4.1}, $i({\cal Y},\pi)=\frac{\underline{j}! a_{\underline{j}}}{\# \text{Aut}({\cal Y})}$, where $\underline{j}=\text{type}(w)$. \hfill${\Box}$

\vskip .5em
As an immediate consequence of Theorem \ref{4.2}, we have:

\begin{Corollary} \label{4.3} Let $\pi, \sigma$ be two projective planes of prime order $p$. Suppose $C_p(\pi)$ and $C_p(\sigma)$ have the same complete weight enumerator. Then, for any partial linear space ${\cal X}$ with at most $\log_2 p$ lines, we have $i({\cal X},\pi)=i({\cal X},\sigma)$.
\end{Corollary}

``{\bf Theorem of Pappus}": Let $\{x_1,x_2,x_3\}$ and $\{y_1,y_2,y_3\}$ be two disjoint 3-sets of collinear points of a projective plane $\pi$. Suppose $\{x_1,x_2,x_3\}$ and $\{y_1,y_2,y_3\}$ determine two distinct lines $\ell_1,\ell_2$ of $\pi$ and the six points $x_i,y_i~ (1\leq i \leq 3)$ are distinct from the point $\ell_1 \wedge \ell_2$. Let us put
\begin{eqnarray*}
z_1 &=& (x_2 \vee y_3) \wedge (x_3 \wedge y_2),\\
z_2 &=& (x_1 \vee y_3) \wedge (x_3 \vee y_1),\\
z_3 &=& (x_1 \vee y_2) \wedge (x_2 \vee y_1).
\end{eqnarray*}
We say that the theorem of Pappus holds in $\pi$, or that $\pi$ is a {\bf Pappian}, if, for every such choice of six initial points $x_i,y_i~ (1 \leq i \leq 3)$ in $\pi$, the three points $z_1,z_2,z_3$ are collinear in $\pi$.

A projective plane need not be Pappian. In fact, a famous theorem in Projective Geometry states (see \cite{6}, \cite{10}) that a projective plane $\pi$ is Pappian iff $\pi$ is the projective plane over a division ring. Since, by the theorem of Wedderburn, the finite division rings are fields, this implies, in particular, that a finite projective plane $\pi$ is Pappian iff it is a field plane. Thus, the finite Pappian planes have prime power orders.

It is easy to see that the nine points $x_i,y_i,z_i ~(1 \leq i \leq 3)$ occuring in Pappus' theorem are necessarily distinct. When Pappus' theorem holds, this set of nine points contains the nine collinear triples listed in Table 1 below. (For some initial choices of the six points $x_i, y_i$, some or all of the three triples $\{x_i, y_i, z_i\}, 1 \leq i \leq 3$, may also be collinear. But this does not affect the following arguments.)

\vskip .5 em

\centerline{{\bf Table 1:} {\bf The lines in Pappus' Configuration}}

\[
\begin{tabular}{c}
  $x_1 x_2 x_3, y_1 y_2 y_3, z_1 z_2 z_3$ \\
  $x_1 y_2 z_3, x_2 y_3 z_1, x_3 y_1 z_2$\\
  $x_1 y_3 z_2, x_2 y_1 z_3, x_3 y_2 z_1$
  \end{tabular}
\]
Consider the partial linear space $\mathbb{P}$ (with nine points and nine lines) which is the subsystem of $PG(2,\mathbb{F}_3)$ obtained as follows. Fix a flag (i.e., an incident point-line pair) $(x,\ell)$ in $PG(2, \mathbb{F}_3)$. Then $\mathbb{P}$  is the subsystem of $PG(2,\mathbb{F}_3)$ whose points are the points of $PG(2,\mathbb{F}_3)$ non-incident with $\ell$, and whose lines are the intersections with this point set of the lines of $PG(2,\mathbb{F}_3)$ non-incident with $x$. Since the automorphism group of $PG(2,\mathbb{F}_3)$ is transitive on the flags, this defines the partial linear space $\mathbb{P}$ uniquely up to isomorphism.

Note that the nine collinear triples of Table 1, occurring in ''the Theorem of Pappus", form an explicit list of the lines of $\mathbb{P}$. This is why $\mathbb{P}$ is sometimes called the {\bf Configuration of Pappus}. (``Configuration" is an old term for a partial linear space.)

Also observe that, despite appearances, the validity (or otherwise) of the ``Theorem of Pappus" does not depend on the explicit ordering of the six initial points, but it depends only on the bijection $x_i \mapsto y_i~ (1 \leq i \leq 3)$ between the initial collinear tuples $\{x_1,x_2,x_3\}$ and $\{y_1,y_2,y_3\}$. More precisely, if the three indices 1, 2, 3 in the statement are consistently permuted, then the validity (or invalidity) of the hypothesis and conclusion of this ``theorem" remains unchanged.

In view of these observations, the theorem of Pappus may be reformulated as follows.

{\bf The theorem of Pappus} (alternative version): Let's say two 3-sets $\alpha, \beta$ of points in a projective plane $\pi$ form an {\bf admissible pair} if (i) $\alpha$ and $\beta$ are collinear triples, (ii) $\alpha$ and $\beta$ are disjoint, and (iii) no four points in $\alpha \sqcup \beta$ are collinear in $\pi$. Then $\pi$ is said to satisfy ``the theorem of Pappus" (or $\pi$ is Pappian) if, for every pair $(\alpha,\beta)$ of admissible triples of $\pi$ and every bijection $f: \alpha \rightarrow \beta$, there is a unique isomorphic copy of $\mathbb{P}$ in $\pi$ such that (a) $\alpha$ and $\beta$ are lines of $\mathbb{P}$, (b) for each $x\in \alpha, x$ and $f(x)$ are non-collinear in $\mathbb{P}$.

Finally, note that the points and lines of $\mathbb{P}$ are uniquely determined by the triple $(\alpha,\beta,f)$ as above. Namely, the nine points and eight of the lines of $\mathbb{P}$ are determined by the hypothesis, and the ninth line of $\mathbb{P}$ is determined by the conclusion of Pappus' Theorem, given $(\alpha, \beta,f)$.

Therefore, the characterization of finite field planes as the finite Pappian projective planes may be rephrased as follows.

\begin{Theorem} \label{4.4} Let $\pi$ be a projective plane of order $n$. Then $i(\mathbb{P},\pi) \leq \frac{2}{3} \binom{n^2+n+1}{2} \binom{n}{3}^2$. Equality holds here iff $\pi$ is a field plane.
\end{Theorem}

{\bf Proof:} Clearly $\pi$ contains exactly $2 \binom{n^2+n+1}{2}\binom{n}{3}^2$ admissible pairs $(\alpha,\beta)$ of collinear triples. For each such pair, there are 3! bijections $f: \alpha \rightarrow \beta$. Thus there are $12 \binom{n^2+n+1}{2} \binom{n}{3}^2$ triples $(\alpha,\beta,f)$ as above. Each such triple determines at most one isomorphic copy of $\mathbb{P}$ in $\pi$. On the other hand, it is easy to verify that $\mathbb{P}$ contains exactly 18 ordered pairs $(\alpha,\beta)$ of disjoint lines, and, for each such pair, there is a unique bijection $f: \alpha \rightarrow \beta$ such that $x$ and $f(x)$ are non-collinear in $\mathbb{P}$ for all $x\in \alpha$. Thus, each subsystem of $\pi$ isomorphic to $\mathbb{P}$ is determined by exactly 18 triples $(\alpha,\beta,f)$. Therefore, a two-way counting argument yields that $i(\mathbb{P},\pi) \leq \frac{12}{18} \binom{n^2+n+1}{2} \cdot \binom{n}{3}^2$, with equality iff $\pi$ is Pappian. \hfill${\Box}$

Using Theorem \ref{4.4} and Corollary \ref{4.3} with ${\cal X}=\mathbb{P}, \sigma =PG(2,\mathbb{F}_p)$, we get (as $\mathbb{P}$ has nine lines):

\begin{Theorem} \label{4.5} Let $\pi$ be a projective plane of prime order $p$ such that $\pi$ has the same complete weight enumerator (of its $p$-ary code) as $PG(2,\mathbb{F}_p)$. If $p > 2^9$, then $\pi$ is isomorphic to $PG(2,\mathbb{F}_p)$.
\end{Theorem}

Recall that a projective plane $\pi$ is said to be {\bf Desarguesian} if (in the standard terminology of projective geometry) each pair of triangles in $\pi$ which is centrally perspective is also axially perspective. Consider the Petersen graph, which may be described as the graph whose vertices are the $\binom{5}{2}$ unordered pairs of symbols from a set of five symbols, with disjointness as adjacency. Let ${\cal D}$ be the partial linear space whose points and lines are both indexed by the vertices of the Petersen graph, such that the line indexed by $x$ is incident with the point indexed by $y$ iff $x$ and $y$ are adjacent vertices of the graph. ${\cal D}$ is known as the {\bf Configuation of Desargue} since it stands in the same relation with ``Desargue's theorem" as $\mathbb{P}$ with the theorem of Pappus". Therefore, the well-known theorem (\cite{6}, \cite{10}) that a projective plane is a field plane iff it is desarguesian may be rephrased as in Theorem \ref{4.4} in the finite case. Namely, for every projective plane $\pi$ of order $n$, one may write down an upper bound for $i({\cal D},\pi)$ in terms of $n$ alone, which is attained iff $\pi$ is a field plane. Using this theorem, one can write an alternative proof of Theorem \ref{4.5}. However, since ${\cal D}$ has ten lines, this alternative proof works only for $p > 2^{10}$. We have chosen to work with $\mathbb{P}$ since it has fewer lines.

\section{Speculations}

The bound $\log_2 p$ in Theorem \ref{4.2} is perhaps the best possible. However, we expect that the bound $p > 2^9$ in Theorem \ref{4.5} is unnecessary, and this theorem actually holds for all primes $p$. For instance, if the conjecture of Neumann (briefly mentioned in the introduction) is correct, then we can use $PG(2,\mathbb{F}_2)$ instead of $\mathbb{P}$ in the proof of Theorem \ref{4.5}, pushing its bound to $p>2^7$. In any case, Theorem \ref{4.5} shows that, in order to prove the uniqueness conjecture \ref{1.2}, at least for large primes $p$, it suffices to calculate the complete weight enumerator for arbitrary projective planes of order $p$. But this is a tall order! We do not even know the complete weight enumerator of $PG(2,\mathbb{F}_p)$ for any prime $p \geq 7$.

To prove that the fourth minimum weight of the $p$-ary code of $PG(2,\mathbb{F}_p)$ is $=3p-3$, we constructed (\cite{3}, \cite{4})~ $p^2(p^3-1) \binom{p+1}{3}$ words of weight $3p-3$ in this code. These words are actually in the dual code. In \cite{2}, \cite{3} we posed

\begin{Conjecture} For any projective plane $\pi$ of prime order $p, C^\perp_p(\pi)$ has at most $p^2(p^3-1)\binom{p+1}{3}$ words of Hamming weight $3p-3$. Equality holds iff $\pi$ is isomorphic to $PG(2,\mathbb{F}_p)$.
\end{Conjecture}

If this conjecture is correct, then, of course, to prove Conjecture \ref{1.2} it will suffice to investigate the initial segment of the Hamming weight enumerator of the dual $p$-ary code of arbitrary projective planes of order $p$.

\end{document}